\documentclass[12pt]{amsart}
\usepackage {amssymb,amsmath}
\makeatletter
\def\@cite#1#2{{\m@th\upshape\bfseries%
[{#1\if@tempswa{\m@th\upshape\mdseries, #2}\fi}]}}
\makeatother
\theoremstyle{plain}
\newtheorem{thm}{Theorem}[section]
\newtheorem{lem}[thm]{Lemma}
\newtheorem{cor}[thm]{Corollary}
\newtheorem{prop}[thm]{Proposition}

\theoremstyle{definition}
\newtheorem{rem}[thm]{Remark}

\newtheorem{note}[thm]{Note}

\newtheorem{eg}[thm]{Example}

\newcommand{\Prf}{\noindent\textbf{Proof.\ }}
\newcommand{\bx}{\strut\hfill$\blacksquare$\medbreak}

\newcommand{\ca}{\mathrm{C}^*}

\newcommand{\bbC}{{\mathbb{C}}}

\newcommand{\bbR}{{\mathbb{R}}}

 \newcommand{\A}{{\mathcal{A}}}
 \newcommand{\B}{{\mathcal{B}}}

 \newcommand{\E}{{\mathcal{E}}}
 
 \newcommand{\G}{{\mathcal{G}}}
\renewcommand{\H}{{\mathcal{H}}}
 
 \newcommand{\J}{{\mathcal{J}}}

 \newcommand{\M}{{\mathcal{M}}}
 
\renewcommand{\O}{{\mathcal{O}}}

 \newcommand{\V}{{\mathcal{V}}}
 \newcommand{\W}{{\mathcal{W}}}

\newcommand{\upchi}{{\raise.35ex\hbox{$\chi$}}}
\newcommand{\qand}{\quad\text{and}\quad}

\newcommand{\qfor}{\quad\text{for}\quad}

\newcommand{\qforal}{\quad\text{for all}\quad}
\newcommand{\qif}{\quad\text{if}\quad}

\newcommand{\Alg}{\operatorname{Alg}}

\newcommand{\rank}{\operatorname{rank}}

\newcommand{\spn}{\operatorname{span}}

\newcommand{\fix}{\operatorname{Fix}}
\def\bra#1{\langle #1|}
\def\ket#1{|#1 \rangle}
\def\one{{\mathchoice{\rm 1\mskip-4mu l}{\rm 1\mskip-4mu l}{\rm 1\mskip-4.5mu l}{\rm
1\mskip-5mu l}}}

\newcommand{\bofh}{\B(\H)}

\begin{document}

\title[Collective Rotation Channels]%
{Noiseless Subsystems for Collective Rotation Channels  in
Quantum Information Theory}
%
\author[J.A. Holbrook,  D.W. Kribs, R. Laflamme,
D. Poulin]{John~A.~Holbrook$^1$,  David~W.~Kribs$^{1,2,3}$,
Raymond Laflamme$^{2,3}$ and David Poulin$^{2,3}$}
\thanks{{\it key words and phrases. } quantum channel, completely positive map, collective rotation channel, quantum error correction, noiseless subsystem, noise commutant.}
\address{$^1$Department of Mathematics and Statistics, University of Guelph,
Guelph, ON, CANADA  N1G 2W1. }
\address{$^2$Institute for Quantum Computing, University of
Waterloo, Waterloo, ON, CANADA N2L 3G1. }
\address{$^3$Perimeter Institute for Theoretical Physics, 35 King
St. North, Waterloo, ON, CANADA N2J 2W9.}

%
\begin{abstract}
Collective rotation channels are a fundamental class of channels
in quantum computing and quantum information theory. The commutant
of the noise operators for such a channel is a $\ca$-algebra which
is equal to the set of fixed points for the channel. Finding the
precise spatial structure of the commutant algebra  for a set of
noise operators associated with a channel is a core problem in
quantum error prevention. We draw on methods of operator algebras,
quantum mechanics and combinatorics to explicitly determine the
structure of the commutant for the class of collective rotation
channels.
\end{abstract}
\maketitle

\section{Introduction}\label{S:intro}

Quantum information theory provides the underlying mathematical
formalism for quantum computing and is an interesting field of
research in its own right \cite{NC}. While quantum computing and
communication promise far reaching applications
\cite{Brooks,Johnson,Nielsen2002}, there are numerous technical
and theoretical difficulties that must be overcome. Of particular
interest is the study of quantum error correction and error
prevention methods. In classical computing, the types of errors
that can occur are very limited. On the other hand, the fragile
nature of quantum systems shows that in quantum computing there is
a much richer variety of potential errors. Fortunately, methods of
quantum error correction have recently been developed showing, in
principle, that these difficulties may be overcome (see
\cite{AB,GottesmanIQEC,Kit,IntroQEC,KLZ,Preskillplenum} for an
introduction to the subject).

Central to quantum information theory is the analysis of {\it
quantum channels} \cite{NC}. Mathematically, a quantum channel is
given by a completely positive trace preserving map which acts on
the set of operators on a finite dimensional Hilbert space. Every
channel has a family of noise operators that determine the map in
a natural way. One of the most promising methods of passive
quantum error correction, recently developed by the third author
and others \cite{DG,Lnoiseless,HKL,IntroQEC,KLV,LCW,ZR}, is called
the {\it noiseless subsystem method}. Given a quantum channel, the
basic tenet  of this method is to use the structure of the
operator algebra defined by the commutant of the associated noise
operators to prepare initial quantum states which are immune to
the noise of the channel. Thus it is a fundamental problem in
quantum error correction to find the structure of this `noise
commutant'. However, let us emphasize that it is the {\it precise}
spatial structure of this algebra that must be identified. This
point is clarified in the discussion of the next section.

An important test class  for the noiseless subsystem method and
other quantum error correction methods is the class of {\it
collective rotation channels}
\cite{BRS,LPS,Fi,Lnoiseless,HKL,KBLW,KLV,VKL,VFPKLC,ZL,Zan,Zar}.
This class has its roots in the depths of quantum mechanics,
specifically in the study of angular momentum at the atomic level
(see for example \cite{QM12}). A realistic physical situation
where these channels arise occurs when quantum information,
encoded as light pulses, is transmitted through an optical fibre
\cite{Lnoiseless,experimental}. In such a situation, the fibre can
produce a `collective rotation' of the information.

In \cite{BS,Kchannel} it was shown that when a channel is unital,
which is the case for collective rotation channels, the noise
commutant is a finite dimensional $\ca$-algebra which is equal to
the fixed point set for the channel. Based on operator algebra
techniques, the paper \cite{HKL} derives an algorithm for
computing the commutant structure in the most general setting.
However, for particular cases such as the channels considered
here, the required computations can become unwieldy.

In this paper, based on the theory of operator algebras and
quantum mechanics, we compute the noise commutant structure for
the class of collective rotation channels. We provide a
constructive proof which yields a simple visual interpretation
based on Pascal's triangle. This result may also be derived from
well-known representation theory techniques; however, our direct
operator theory cum quantum mechanical approach is novel and
offers a new perspective on the general problem.

The next section contains a brief review of the material we
require from the theories of operator algebras and quantum
information. In the third section we define the collective
rotation channels and establish some basic properties. The fourth
section contains the commutant structure theorem for the `qubit'
case (Theorem~\ref{commthm}).  Finally, we conclude the paper by
presenting a commutant structure theorem for more general classes
of collective rotation channels (Theorem~\ref{commthm1}).

One final comment. A study of the quantum information and quantum
computing  literature reveals that many techniques from operator
theory and operator algebras have been, or could be, used to build
mathematical foundations for the physical theories in these areas.
An idea we wish to promote is that there is a wealth of
interesting mathematics to be found in this young field.

\section{Background}\label{S:back}

Motivated by the postulates of quantum mechanics, an assumption
typically made in quantum information theory is that every quantum
operation on a closed quantum system is reversible \cite{QM12,NC}.
Mathematically, this statement means that the operation is
described by  unitary evolution; in other words, there is a
unitary operator $U$ on a Hilbert space $\H$ such that the
operation is implemented by the conjugation map $\rho \mapsto U
\rho U^\dagger$ where $\rho$ is an operator  on $\H$. (Here we use
the physics convention $U^\dagger$ for conjugate transpose.) Often
$\rho$ is a density operator, a positive operator with trace equal
to one, that corresponds to the initial state of the quantum
system of interest, but in our analysis there is no loss of
generality in considering evolution of any operator under the
quantum operation. Further note that $U$ can be restricted to the
special unitary group $SU(N)$, where $N = \dim(\H)$, since the
evolution $\rho \mapsto U \rho U^\dagger$ is unaffected by the
multiplication of $U$ by a complex phase.

Of course, in practice a given quantum operation will not be
reversible because of interactions with the environment. In this
more realistic setting the quantum operation is regarded as acting
on a closed quantum system that contains the original as a
subsystem. The mathematical formalism for this is given by
completely positive maps
\cite{Choi,Kraus,Paulsentext,Paulsentext2} and the Stinespring
dilation theorem \cite{Stine}. Specifically, every quantum
operation is represented mathematically by a quantum channel.

Given a (finite dimensional) Hilbert space $\H$, a {\it quantum
channel} is a map $\E$ which acts on the set $\bofh$ of all
operators on $\H$ and is completely positive and trace preserving.
For each channel $\E$ there is a set of (non-unique) {\it noise
operators} \cite{Choi,Kraus} $\{ A_1, \ldots,A_n\}$ that determine
the map through the equation
\begin{eqnarray}\label{channeldef}
\E(\rho) = \sum_{k=1}^n A_k \rho A_k^\dagger \qfor \rho \in\bofh.
\end{eqnarray}
Physically, the associated quantum operation can be regarded as
determined by a compression of the Stinespring unitary dilation,
that acts on a larger closed quantum system, of the completely
positive map (\ref{channeldef}). Trace preservation is equivalent
to the noise operators satisfying the equation
\[
\sum_{k=1}^n A_k^\dagger A_k = \one,
\]
where $\one$ is the identity operator on $\H$. The channel is
unital if also,
\[
\E(\one) = \sum_{k=1}^n A_k A_k^\dagger  = \one.
\]

Let $\fix(\E) = \{ \rho \in \bofh : \E(\rho) = \rho \}$ be the
fixed point set for $\E$ and let $\A$ be the algebra generated by
$A_1, \ldots, A_n$ from (\ref{channeldef}). This is called the
{\it interaction algebra} in quantum information theory
\cite{KLV}. In general, $\fix(\E)$ is just a $\dagger$-closed
subspace of $\bofh$, but it was shown (independently) in \cite{BS}
and \cite{Kchannel} that, in the case of a unital channel $\E$,
the so-called {\it noise commutant} $\A^\prime = \{ \rho \in \bofh
: \rho A_k =A_k \rho, \,\, k=1,\ldots,n\}$ coincides with this
set:
\[
\fix(\E) = \A^\prime.
\]
In particular, $\fix(\E) = \A^\prime$ is a $\dagger$-closed
operator algebra (a finite dimensional $\ca$-algebra
\cite{Arvinvite,byeg}). Further, the von Neumann double commutant
theorem shows how the algebra $\A = \A^{\prime\prime} =
\fix(\E)^\prime$ only depends on the channel; that is, it is
independent of the choice of noise operators that determine the
channel as in (\ref{channeldef}).

It is a fundamental result in finite dimensional $\ca$-algebra
theory \cite{Arvinvite,byeg,Tak} that every such algebra is
unitarily equivalent to an orthogonal direct sum of `ampliated'
full matrix algebras; i.e., there is a unitary operator $U$ such
that
\[
U \A U^\dagger = \sum_{k=1}^d \oplus \,\,\big( \one_{m_k}\otimes
\M_{n_k}  \big),
\]
where $\M_{n_k}$ is the full matrix operator algebra
$\B(\bbC^{n_k})$. The numbers $m_k$ in this decomposition
correspond to the multiplicities in the $\ca$-algebra
representation that gives $\A$.  With this form for $\A$ given,
the structure of the commutant up to unitary equivalence is easily
computed by
\begin{eqnarray}\label{commform}
\A^\prime \simeq \sum_{k=1}^d \oplus\,\, \big( \M_{m_k} \otimes
\one_{n_k} \big).
\end{eqnarray}
(See \cite{Lnoiseless,HKL,KBLW,VKL,ZL,Zan} for  more detailed
discussions in connection with quantum information theory.)

On the other hand, for a given quantum channel $\E$ with noise
operators $\{A_k\}$, the noise commutant $\A^\prime$ plays a
significant role in quantum error prevention. The structure of
this commutant can be used to prepare density operators, which
encode the state of a given quantum system, for use in the
noiseless subsystem method of error correction. This is a passive
method of quantum error correction, in the sense that such
operators will remain immune to the effects of the noise
operators, or `errors' of the channel, without active
intervention. But more is true. The algebra structure discussed
above shows that quantum operations may be performed on such a
subsystem,  provided the corresponding unitary operators belong to
the commutant. Keeping in mind our earlier description of an
optical fibre, the reader can imagine a situation where it is
desirable to transfer quantum information through the fibre such
that the information remains immune to the errors of collective
rotations produced by the fibre.

As discussed above, understanding the structure of $\A^\prime$ is
of fundamental importance in quantum error correction. But there
is an operator algebra subtlety here which is worth emphasizing.
Typically, it is not feasible in this setting to wash away the
particular representation which gives $\A^\prime$ with
$\ast$-isomorphisms, unitary equivalences, etc., as is the custom
in operator algebra theory. Indeed, by the very nature of the
problems, it is the precise spatial algebra structure of
$\A^\prime$ which must be identified, ampliations included.

The basic problem of computing $\A^\prime$ was addressed in
\cite{HKL} for the general case of a unital quantum channel. We
also mention more recent work \cite{Zar} where computer algorithms
have been written for this and other related purposes. However, in
particular cases, such as the class of channels considered in this
paper, a more delicate approach based on special properties of the
class can be exploited to find this structure more directly and
efficiently.

\section{Collective Rotation Channels}\label{S:coll}

Let $\{ \ket{\frac{-1}{2}}, \ket{\frac{1}{2}}\}$ be a fixed
orthonormal basis for 2-dimensional Hilbert space $\H_2 = \bbC^2$,
corresponding to the classical base states in a two level quantum
system (e.g. the ground and excited states of an electron in a
Hydrogen atom). Note that such a basis is usually written as $\{
\ket{0}, \ket{1}\}$, but the $\frac {-1}2 , \frac{1}2$ notation is
more convenient for the combinatorics below. A `qubit' or `quantum
bit' of information is given by a unit vector $\ket{\psi} = \alpha
\ket{\frac {-1}2} + \beta \ket{\frac{1}2}$ inside $\H_2$. When
both $\alpha$ and $\beta$ are non-zero, $\ket{\psi}$ is said to be
a {\it superposition} of $\ket{\frac {-1}2}$ and
$\ket{\frac{1}2}$.

We shall make use of the abbreviated form from quantum mechanics
for the associated standard orthonormal basis for $\H_{2^n} =
(\bbC^2)^{\otimes n} \simeq \bbC^{2^n}$. For instance, the basis
for $\H_4$ is given by
\[
\Big\{ \ket{ij} : i,j\in\{ \frac {-1}2 , \frac{1}2 \} \Big\}
\]
where $\ket{ij}$ is the vector tensor product $\ket{ij} \equiv
\ket{i}\ket{j} \equiv \ket{i}\otimes\ket{j}$.

Let $\{ \sigma_x, \sigma_y, \sigma_z \}$ be the spin-$1/2$ Pauli
matrices given by
\[
\sigma_x = \frac{1}{2} \left(\begin{matrix}
0&1 \\
1&0
\end{matrix}\right), \quad
\sigma_y = \frac{1}{2} \left(\begin{matrix}
0&-i \\
i&0
\end{matrix}\right), \quad
\sigma_z = \frac{1}{2} \left(\begin{matrix}
1&0 \\
0&-1
\end{matrix}\right).
\]
Further let $\one_2$ be the $2\times 2$ identity matrix. We shall
regard these as the matrix representations for operators acting on
$\H_2$ with respect to $\{ \ket{\frac {-1}2}, \ket{\frac{1}2}\}$.
The Pauli matrices satisfy the following commutation relations:
\begin{eqnarray*}
(1) & [\sigma_x,\sigma_y] = i \sigma_z \\
(2) & [\sigma_z,\sigma_x] = i\sigma_y \\
(3) & [\sigma_y,\sigma_z] = i\sigma_x.
\end{eqnarray*}

These are the canonical commutation relations which define the Lie
algebra $su(2)$, given by the linear space $r_x \sigma_x + r_y
\sigma_y + r_z \sigma_z = \vec r \cdot \vec \sigma$ with
$(r_x,r_y,r_z) \in \bbR $.  This algebra is the generator of the
Lie group $SU(2)$ as the manifold of $2\times 2$ unitary matrices
with unit determinant and is isomorphic to the manifold $\{\exp(-i
2\pi\vec r \cdot \vec \sigma) : ||\vec{r}|| \leq 1\}$. The group
$SU(2)$ is referred to as the rotation group as it is homeomorphic
to $O(3)$, the rotational group in three-dimensional space. Note
that a rotation is the most general transformation which can be
performed on a closed two-dimensional quantum system.

Now let $n\geq 1$ be a fixed positive integer. Define operators
$\{ J_z^{(k)} : 1\leq k \leq n\}$ on $\H_{2^n}$ by
\[
J_z^{(1)} = \sigma_z \otimes (\one_2)^{\otimes (n-1)}, \quad
J_z^{(2)} = \one_2 \otimes \sigma_z \otimes (\one_2)^{\otimes
(n-2)}, \,\,\,\,\,\ldots,
\]
where we use  the standard ordering $(a_{kl}B)_{kl}$ for the
tensor product of  matrices $A\otimes B$. Similarly define $\{
J_x^{(k)}, J_y^{(k)} : 1\leq k \leq n\}$. Then the {\it collective
rotation  generators} $\{ J_x, J_y, J_z\}$ are given by
\[
J_x = \sum_{k=1}^n J_x^{(k)}, \quad J_y = \sum_{k=1}^n J_y^{(k)},
\quad J_z = \sum_{k=1}^n J_z^{(k)}.
\]
Let us set down the fundamental commutation relations satisfied by
these operators \cite{QM12}.

\begin{prop}\label{liealg}
The following relations hold for $\{J_x, J_y, J_z\}$:
\begin{eqnarray*}
(1)^\prime & [J_x,J_y] = iJ_z \\
(2)^\prime & [J_z,J_x] = iJ_y \\
(3)^\prime & [J_y,J_z] = iJ_x.
\end{eqnarray*}
\end{prop}

\Prf These identities easily follow from corresponding equations
for $J_x^{(k)},J_y^{(k)},J_z^{(k)}$, with $1\leq k\leq n$, which
are  simple consequences of the commutation relations
$(1),(2),(3)$ satisfied by the spin-$1/2$ Pauli matrices. \bx

\begin{note}\label{repn}
Observe that Proposition~\ref{liealg} shows $\{J_x, J_y, J_z\}$
determine a $2^n$-dimensional representation of $su(2)$.
\end{note}

In what follows, much of the analysis will be focused on the
operators
\[
J_+ = J_x + iJ_y \qand J_- = J_x - iJ_y =J_+^\dagger.
\]
We shall also consider the so-called {\it $J$-total} operator
$J^2$ defined by
\[
J^2 = J_x^2 + J_y^2 + J_z^2.
\]
The $J^2$ notation comes from the fact that this operator is
conventionally  defined as a vector product of matrices
\cite{QM12}.

Intuitively, the collective rotation channel is one where every
qubit undergoes the same unknown rotation. Let us formalize this
notion. Consider a channel $\E_{U^{\otimes n}}: \B(\H_{2^n})
\rightarrow \B(\H_{2^n})$ defined as $\E_{U^{\otimes n}}(T) =
\left( U^{\otimes n}\right) T \left(U^{\otimes n}\right)^\dagger$
for $U \in SU(2)$. This is a collective rotation of $n$ qubits
which can also be written $\E_{U^{\otimes n}}(T) = \exp(-i
2\pi\vec r \cdot \vec J) T \exp(i 2\pi\vec r \cdot \vec J)$, where
$\vec r \cdot \vec J = r_x J_x + r_y J_y + r_z J_z$ and $U = \exp
(-i 2 \pi \vec{r} \cdot \vec{\sigma})$. Hence the appellation {\em
collective rotation generators}.

But here, the specific rotation $U$ is unknown and chosen at
random  over $SU(2)$ according to some probability distribution
$P(\vec{r})$, for instance the distribution corresponding to Haar
measure on $SU(2)$. Hence, the {\em $n$-qubit collective rotation
channel} can be written as
\begin{equation}
\E_n(T) = \int_{\{||\vec{r}||\leq 1\}} \exp(-i 2\pi\vec r \cdot
\vec J) \,\,T\,\, \exp(i 2\pi\vec r \cdot \vec J) P(\vec{r}) d
\vec r;
\end{equation}
it is a weighted average of all collective rotations. By the
symmetry of the integrated region, it can be shown that this
unital channel can also be expressed in a more conventional form,
\begin{eqnarray}\label{nqubitchannel}
\E_n (T) = E_x T E_x^\dagger + E_y T E_y^\dagger + E_z T
E_z^\dagger,
\end{eqnarray}
where the noise operators are defined as
\[
E_x =\frac{1}{\sqrt{3}}\exp (i\theta_x J_x),\quad E_y=
\frac{1}{\sqrt{3}}\exp (i\theta_y J_y), \quad E_z=
\frac{1}{\sqrt{3}}\exp (i\theta_z J_z),
\]
and $\theta_k$, $k=x,y,z$, are angles determined by the
probability distribution.


It is not hard to see that our analysis is independent of the
particular choices for these angles, provided each $\theta_k$ is
non-zero. Indeed, through a standard functional calculus argument
from operator algebra, it can be seen that the interaction
algebras generated by the $J_k$ and $E_k$ coincide, whatever the
choice of $\theta_k$;
\begin{eqnarray}\label{Adefn}
\A_n \equiv \Alg \{J_x, J_y, J_z\} =\Alg\{J_+,J_-,J_z\} = \Alg
\{E_x, E_y, E_z\}.
\end{eqnarray}
In particular, as observed in \cite{HKL}, the fixed point set of
this channel is determined by the original rotation generators.

\begin{prop}
Let $n\geq 1$ be a positive integer. Then
\begin{eqnarray*}
\fix (\E_n) =\A_n^\prime = \{E_x, E_y, E_z \}^\prime = \{J_x, J_y,
J_z \}^\prime.
\end{eqnarray*}
Further, this commutant may be computed by considering the joint
commutant of any pair from $\{J_x, J_y, J_z \}$.
\end{prop}

\section{Commutant Structure Theorem}\label{S:comm}

Given a positive integer $n\geq 1$, let  $\Delta_n$ denote  the
graph of ${n \choose \cdot}$; that is, the graph of the $n$th line
in Pascal's triangle. (See the example below for a pictorial
perspective.) Let
\[
\J_n = \left\{
\begin{array}{cl}
\{0,1,\ldots, \frac{n}{2} \} & \mbox{if $n$ is even} \\
\{\frac{1}{2}, \frac{3}{2}, \ldots, \frac{n}{2} \} & \mbox{if $n$
is odd}
\end{array}\right.
\]
Observe that the cardinality of $\J_n$ is equal to the number of
steps up one side of $\Delta_n$.

\begin{thm}\label{commthm}
Let $\E_n$ be the collective rotation channel for a fixed positive
integer $n\geq 1$. Then
\begin{eqnarray}
\fix (\E_n) = \A_n^\prime = \sum_{j\in\J_n} \oplus \,
\A_{(j)}^\prime,
\end{eqnarray}
where $\A_{(j)}^\prime$ is a $\ca$-subalgebra of $\A_n^\prime$
given, up to unitary equivalence, by
\[
\A_{(j)}^\prime \simeq \M_{p_j} \otimes \one_{q_j} \qfor j\in\J_n,
\]
with $p_{\frac{n}{2}} = 1$ and for $j\in\J_n$, $j< \frac{n}{2}$,
\[
p_j = {n \choose j + \frac{n}{2}} - {n \choose j +\frac{n}{2} + 1}
= {n+1 \choose j + \frac{n}{2} + 1} \frac{q_j}{n+1},
\]
where
\[
q_j = 2j +1 \qfor j\in\J_n.
\]
\end{thm}

In the proof below we shall explicitly identify the spatial
decomposition that yields this decomposition of $\A^\prime_n$.
Recall that this is necessary for using the noiseless subsystem
approach to quantum error correction. Before proving this theorem,
let us illustrate how $\Delta_n$ gives a visual method for
determining the commutant structure. For the sake of brevity, let
us focus on a single case, the $n=4$ collective rotation channel
$\E_4$.

\begin{eg}\label{4qubiteg}
In the $n=4$ case we have $\J_4 = \{ 0,1,2\}$ and $p_0=2, p_1=3,
p_2=1$ and $q_0=1, q_1=3, q_2=5$. The theorem states that
\[
\fix (\E_4) = \A_4^\prime = \A_{(0)}^\prime \oplus \A_{(1)}^\prime
\oplus \A_{(2)}^\prime,
\]
with each $\A^{(j)}$ a subalgebra of $\A_4^\prime$ unitarily
equivalent to
\begin{eqnarray*}
\A_{(0)}^\prime &\simeq& \bbC \otimes \one_5 \simeq \bbC \one_5 \\
\A_{(1)}^\prime &\simeq& \M_3 \otimes \one_3  \\
\A_{(2)}^\prime &\simeq& \M_2 \otimes \one_1 \simeq \M_2.
\end{eqnarray*}
Consider the structure of $\Delta_4$:

\vspace{0.1in}

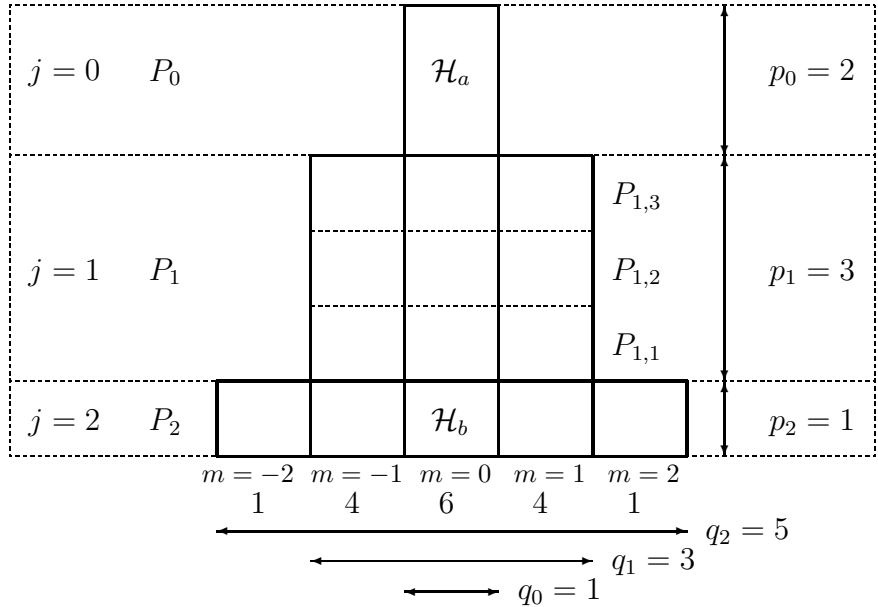
\setlength{\unitlength}{0.05cm}
\begin{picture}(230,160)

{ \thicklines

\put(105,40){\line(0,1){120}}

\put(130,40){\line(0,1){120}}

\put(105,160){\line(1,0){25}}

\put(80,40){\line(0,1){80}}

\put(155,40){\line(0,1){80}}

\put(80,120){\line(1,0){75}}

\put(55,40){\line(0,1){20}}

\put(180,40){\line(0,1){20}}

\put(55,60){\line(1,0){125}}

\put(55,40){\line(1,0){125}}

}

\put(112,140){$\H_a$}

\put(112,47){$\H_b$}



\put(80,80){\dashbox{1}(75,20)}

\put(0,40){\dashbox{1}(230,20)}

\put(0,60){\dashbox{1}(230,60)}

\put(0,120){\dashbox{1}(230,40)}

\put(160,67){$P_{1,1}$}

\put(160,87){$P_{1,2}$}

\put(160,107){$P_{1,3}$}

\put(37,140){$P_0$}

\put(37,87){$P_1$}

\put(37,47){$P_2$}

\put(5,140){$j=0$}

\put(5,87){$j=1$}

\put(5,47){$j =2$}

\put(202,140){$p_0=2$}

\put(202,87){$p_1=3$}

\put(202,47){$p_2 =1$}

\put(190,40){\vector(0,1){20}}

\put(190,60){\vector(0,-1){20}}

\put(190,60){\vector(0,1){60}}

\put(190,120){\vector(0,-1){60}}

\put(190,120){\vector(0,1){40}}

\put(190,160){\vector(0,-1){40}}

\put(51,33){{\Small $m=-2$}}

\put(64,25){$1$}

\put(80,33){{\Small $m=-1$}}

\put(89,25){$4$}

\put(109,33){{\Small $m=0$}}

\put(114,25){$6$}

\put(134,33){{\Small $m=1$}}

\put(139,25){$4$}

\put(159,33){{\Small $m=2$}}

\put(164,25){$1$}

\put(55,20){\vector(1,0){125}}

\put(180,20){\vector(-1,0){125}}

\put(185,18){$q_2=5$}

\put(160,10){$q_1=3$}

\put(80,12){\vector(1,0){75}}

\put(155,12){\vector(-1,0){75}}

\put(105,4){\vector(1,0){25}}

\put(130,4){\vector(-1,0){25}}

\put(135,2){$q_0=1$}

\end{picture}

The number $p_j$ corresponds to the `height' of the $j$th
horizontal bar (counting top-down), and $q_j$ equals the number of
blocks inside this bar. Spatially, the vertical bars correspond to
the eigenspaces for $J_z$ for the eigenvalues $m=-2,-1,0,1,2$
(with eigenspace projections $Q_m$ in the proof below), which have
respective multiplicities 1, 4, 6, 4, 1. The horizontal bars
correspond to eigenspaces of $J^2$ (Corollary~4.10). The
corresponding eigenspace projections $P_0,P_1,P_2$ are the minimal
central projections for $\A_4$ and $\A_4^\prime$.

To see how the blocks correspond to subspaces, the subspace
$\H_a$, as an example, for the top box in $\Delta_4$ is the joint
eigenspace for $J_z$ and $J^2$, corresponding to $m=0$ and $j=0$
with our notation below. Each of the $j$th horizontal bars further
breaks up into smaller horizontal bars, for instance $P_1 =
\sum_{k=1}^3 P_{1,k}$. The subspaces $\{P_{j,k}\H\}$ form the
maximal family of minimal reducing subspaces for $\A_4$ as
outlined below. On the other hand, the corresponding family for
$\A_4^\prime$ is given by the vertical blocks inside the $j$th
horizontal bar. For example, the projection onto $\H_b$ and the
projections onto its other four counterparts in the $j=2$ bar
(which are all 1-dimensional because they lie in the $j=2$ bar)
are the family of minimal $\A_4^\prime$-reducing subspaces
supported on $P_2$.
\end{eg}

We now turn to the proof of Theorem~\ref{commthm}. Let $n\geq 1$
be a fixed positive integer. We shall find the structure of
$\A^\prime_n$ by first computing the structure of $\A_n$. We begin
by showing how the numeric distribution of the eigenvalues for
$J_z$ is linked with $\Delta_n$. In what follows, we use the
abbreviated Dirac notation to denote the standard orthonormal
basis for $\H \equiv \H_{2^n} = \bbC^{2^n}$ with
$\ket{\frac{-1}{2}}, \ket{\frac{1}{2}}$ corresponding to the base
states of the two-level quantum system ($d=2$ with our notation in
the next section);
\[
\Big\{ \ket{\vec{i}} = \ket{i_1 i_2 \cdots i_n} : i_j \in
\{\frac{-1}{2}, \frac{1}{2}\},\, 1\leq j \leq n\Big\}.
\]

\begin{lem}\label{Jzeigen}
For $m = -\frac{n}{2}, -\frac{n}{2}+1,\ldots, \frac{n}{2}$
consider the subspaces of $\H$ given by
\[
\V_m = \spn \big\{ \ket{\vec{i}} : | \vec{i} | =m  \big\},
\]
where $|\vec{i}| = \sum_{j=1}^n i_j$. Then $\H =
\sum_{m=-\frac{n}{2}}^{\frac{n}{2}} \oplus \, \V_m$ and
\[
\dim \V_m = {n \choose m + \frac{n}{2}} \qfor -\frac{n}{2}\leq m
\leq \frac{n}{2}.
\]
Further, $\V_m$ is an eigenspace for $J_z$ corresponding to the
eigenvalue
\[
\lambda = m  \qfor -\frac{n}{2} \leq m \leq \frac{n}{2}.
\]
\end{lem}

\Prf The spatial decomposition of $\H$ is easy to see and the
dimensions of the $\V_m$ follow from simple combinatorics. For the
eigenvalue connection with $J_z$, observe that for $|\vec{i}| =m$
we have
\[
J_z \ket{\vec{i}} = \sum_{k=1}^n J_z^{(k)} \ket{\vec{i}} =
\sum_{k=1}^n i_k \ket{\vec{i}}  = |\vec{i}| \ket{\vec{i}} = m
\ket{\vec{i}}.
\]
\bx

For $-\frac{n}{2}\leq m \leq \frac{n}{2}$, let $Q_m$ be the
orthogonal projection of $\H$ onto $\V_m \equiv Q_m \H$.

\begin{lem}\label{jplusjz}
Given $-\frac{n}{2}\leq m \leq \frac{n}{2}$, we have
\[
J_+ Q_m = \left\{ \begin{array}{cl}
Q_{m+1} J_+ Q_m & \mbox{if $m < \frac{n}{2}$} \\
0 & \mbox{if $m = \frac{n}{2}$} \end{array}\right.
\]
and
\[
J_- Q_m = \left\{ \begin{array}{cl}
Q_{m-1} J_- Q_m & \mbox{if $m > -\frac{n}{2}$} \\
0 & \mbox{if $m = - \frac{n}{2}$}
\end{array}\right.
\]
\end{lem}

\Prf Let $\ket{\psi}$ belong to $Q_m\H$. Then $J_z \ket{\psi} = m
\ket{\psi}$. But notice that
\begin{eqnarray*}
J_z J_+ = J_z (J_x + iJ_y) &=& J_xJ_z +iJ_y + iJ_yJ_z +J_x \\ &=&
J_+ (J_z + \one).
\end{eqnarray*}
Thus $J_z J_+ \ket{\psi} = (m +1) J_+ \ket{\psi}$ when $m<
\frac{n}{2}$, so that $J_+ \ket{\psi}$ belongs to $Q_{m+1}\H$. The
corresponding identities for $J_-$ are proved in a similar fashion
and for convenience the identities $J_+ Q_m = 0 = J_- Q_{-m}$, $m
= \frac{n}{2}$, will be observed in the discussion which follows.
\bx

Next we shall derive a spatial decomposition of $\H$ which will
allow us to connect with the structure of $\Delta_n$. Let
\[
\ket{0_L}\equiv \ket{\frac{n}{2}, - \frac{n}{2},1}
\]
be a  (unit) eigenvector for $J_z$ for the eigenvalue $m =
-\frac{n}{2}$. The span of $\ket{0_L}$ will be identified with the
`bottom left corner' of $\Delta_n$, see
Corollary~\ref{spatialdecomp} below. To simplify notation, let $ns
= \frac{n}{2}$ (the use of this notation will become clear in the
next section). Lemma~\ref{jplusjz} shows that $J_+ \ket{0_L}
\equiv \ket{ns, -ns +1,1}$ is an eigenvector of $J_z$ for the
eigenvalue $m= - ns + 1$. Similarly, the vectors
\[
J_+^p \ket{0_L} \equiv \ket{ns, -ns + p, 1} \qfor 0\leq p <
q_{ns},
\]
are non-zero and belong to $\V_{-ns+ p}$.

Let $\{ \ket{ns-1,-ns+1,\mu}\}_\mu$ be an orthonormal basis for
$\V_{-ns+1} \ominus \spn\{ \ket{ns,-ns+1,1}\}$.  Now inductively,
if we are given $j\in\J_n$ with $j < ns$, let $\ket{j,m=-j,\mu}$
be an orthonormal basis for
\[
\V_{-j} \ominus \spn\{ \ket{j^\prime, -j,\mu}: j< j^\prime \leq
ns\},
\]
where $\ket{j^\prime, -j,\mu} = J_+^{(j^\prime - j)}
\ket{j^\prime, -j^\prime,\mu}.$

Notice that
\begin{eqnarray}\label{jminus}
J_- \ket{j,m=-j,\mu} = 0 \qforal j,\mu.
\end{eqnarray}
Indeed, by choice of the vectors $\ket{j,-j,\mu}$ and from the
`eigenspace shifting' of Lemma~\ref{jplusjz}, it follows that each
$\ket{j,-j,\mu}$ is orthogonal to the range space of $J_+$. Thus
$J_-$ annihilates the left hand steps of $\Delta_n$, which is the
content of (\ref{jminus}). From this we also have
\begin{eqnarray}\label{jplus}
J_+  \ket{j,m=j,\mu} = J_+^{2j} \ket{j,m=-j,\mu}= 0 \qforal \mu.
\end{eqnarray}
In other words, from the $\Delta_n$ picture given by
Corollary~\ref{spatialdecomp} below, $J_+$ annihilates the right
hand side blocks of $\Delta_n$.

Thus, in summary we have a collection of vectors $\ket{j,m,\mu}$
(which turn out to form an orthogonal basis for $\H$) such that:
$j$ belongs to $\J_n$, for fixed $j$ the range of $m$ is $-j\leq m
\leq j$, for each $j,m,\mu$,
\[
\ket{j,m,\mu} = J_+^{(m+j)} \ket{j,-j,\mu} = J_-^{(m-j)}
\ket{j,j,\mu},
\]
and for a given $m,j$ pair the index $\mu$ has $p_j$ possible
values.

For fixed $j,\mu$ let $\H(j,\mu)$ be the subspace defined by
\[
\H(j,\mu) = \spn\big\{ \ket{j,m,\mu}  : -j\leq m \leq j\big\}.
\]
Such a subspace corresponds to a horizontal slice of the `$j$th
horizontal bar' in $\Delta_n$. From Corollary~\ref{spatialdecomp},
it follows that these subspaces are pairwise orthogonal for
distinct pairs $j,\mu$. (This justifies the use of the orthogonal
sum symbol $\oplus$ in the following statement.)

\begin{lem}\label{J2eigen}
The operator $J^2$ belongs to the centre of $\A$; that is,
\[
J^2 \in\A_n\cap \A_n^\prime.
\]
Consider the subspaces
\[
\W_j = \sum_\mu \oplus\,\, \H(j,\mu) \qfor j\in\J_n.
\]
Then the restriction of $J^2$ to each of  these subspaces is a
constant operator; i.e., there are scalars $\lambda_j$ such that
\[
J^2|_{\W_j} = \lambda_j\one_{\W_j} \qfor j\in\J_n.
\]
Further, these scalars satisfy $\lambda_{j_1}\neq \lambda_{j_2}$
for $j_1 \neq j_2$.
\end{lem}

\Prf By definition $J^2$ belongs to $\A$. We show that $J^2$
commutes with $J_x$. The $J_y$ and $J_z$ cases are similar.
Observe that
\begin{eqnarray*}
[J_x,J^2] &=& [J_x, J^2_y + J^2_z] =[J_x, J_y^2] + [ J_x,J_z^2]
\\
&=& J_y[J_x,J_y] + [J_x, J_y]J_y + J_z[J_x,J_z] + [J_x, J_z]J_z \\
&=& J_y (iJ_z) + (iJ_z) J_y + J_z (-iJ_y) + (-i J_y) J_z = 0.
\end{eqnarray*}

Consider a  vector $\ket{j,-j,\mu}$ in the left most block of
$\W_j$. Observe that $ J^2 = J_+J_- + J_z^2 -J_z, $ and hence
\[
J^2 \ket{j,-j,\mu} = (J^2_z -J_z)\ket{j,-j,\mu}  = (j^2+j)
\ket{j,-j,\mu}.
\]
As $J^2$ belongs to $\A_n^\prime$, we have $J^2J_+ = J_+J^2$.
Thus, given a typical basis vector $J_+^{(m+j)}\ket{j,-j,\mu} =
\ket{j,m,\mu}$ inside $\W_j$ compute
\begin{eqnarray*}
J^2(J_+^{(m+j)}\ket{j,-j,\mu}) &=& J^{(m+j)}_+ J^2 \ket{j,-j,\mu}
\\ &=& (j^2+j)J_+^{(m+j)}\ket{j,-j,\mu} = (j^2+j)\ket{j,m,\mu}.
\end{eqnarray*}
It follows that the corresponding restrictions of $J^2$ satisfy
$J^2|_{\W_j} = (j^2+j) \one_{\W_j}$, and the scalars $\lambda_j =
j^2+j$ are different for distinct values of $j$. \bx


\begin{cor}\label{spatialdecomp}
The vectors $\{\ket{j,m,\mu}  \}_{j,m,\mu}$ are non-zero and form
an orthogonal basis for $\H$. Thus,
\[
\H = \sum_j \oplus \,\W_j = \sum_{j,\mu} \oplus \,\H(j,\mu),
\]
and the subspaces $\{ \W_j \}$ are the eigenspaces for $J^2$.
\end{cor}

\Prf These vectors are clearly all non-zero by the above
discussions. Consider two vectors from this set,
$J_+^{p_i}\ket{j_i,m_i,\mu_i}$ for $i=1,2$. Then
\[
\bra{j_1,m_1,k_1} J_+^{p_2-p_1} \ket{j_2,m_2,k_2} = 0 \qif
(j_1,m_1,k_1) \neq (j_2,m_2,k_2).
\]
This follows from the choice of the vectors $\ket{j,m,\mu}$, the
relations
\begin{eqnarray}\label{jplusjminus}
J_+ J_- =  J_x^2 + J_y^2 - J_z = J^2 - J_z^2 -J_z,
\end{eqnarray}
\begin{eqnarray}\label{jminusjplus}
J_- J_+ =  J_x^2 + J_y^2 + J_z = J^2 - J_z^2 +J_z,
\end{eqnarray}
and the connections with the eigenspaces for $J_z, J^2$ given by
Lemma~\ref{Jzeigen} and Lemma~\ref{J2eigen}. \bx

The following perspective on the actions of $J_+$ and $J_-$ will
be useful below.

\begin{lem}\label{weighted}
For all $j,\mu$, the operators $J_+$ and $J_- = J_+^\dagger$ act
as weighted shifts on the standard basis for $\H(j,\mu)$.
\end{lem}

\Prf Recall that $J_- \ket{j,m=-j,\mu} = 0$ since $\ket{j,-j,\mu}$
belongs to the orthocomplement of the range of $J_-^\dagger =
J_+$; that is, $\bra{j,-j,\mu}\ket{J_+\psi} =0$ for all
$\ket{\psi} \in \H$. Thus, by equation (\ref{jminusjplus}) and
Lemma~\ref{Jzeigen} and Lemma~\ref{J2eigen}, for $p\geq 1$ there
is a scalar $c$ with
\begin{eqnarray*}
J_- J_+^{(m+j)} \ket{j,m,\mu} &=& (J^2 -J_z^2 -J_z) J_+^{(m+j-1)} \ket{j,m,\mu} \\
&=& c J_+^{p-1} \ket{j,m-1,\mu}.
\end{eqnarray*}
In particular, $J_-$ acts as a backward shift on the (orthogonal)
basis $\{ \ket{j,m,\mu} : -j\leq m \leq j\}$ for $\H(j,\mu)$ with
$J_- \ket{j,-j,\mu} = 0$. Similarly, by using (\ref{jplusjminus})
it can be seen that $J_+$ acts as the forward shift on this basis
with $J_+ \ket{j,j,\mu} =0$.

Hence, when the basis $\{\ket{j,m,\mu} : -j \leq m \leq j\}$ is
normalized to turn it into an orthonormal basis for $\H(j,\mu)$,
we see that $J_+$ (respectively $J_-$) acts as a forward
(respectively backward) weighted shift on this basis. \bx

The following result shows that the family of mutually orthogonal
subspaces $\H(j,\mu)$ forms the (unique) maximal family of minimal
reducing subspaces for $\A_n$ which determine the minimal central
projections.

\begin{lem}\label{minreduce}
For all $j,\mu$, the subspace $\H(j,\mu)$ is a minimal
$\A_n$-reducing subspace.
\end{lem}

\Prf First note that $\H(j,\mu)$ is clearly reducing for $J_z$
(i.e. invariant for both $J_z$ and $J_z^\dagger$). Also,
Lemma~\ref{weighted} shows that $\H(j,\mu)$ reduces $J_+$ and
$J_-$. Hence $\H(j,\mu)$ is a reducing subspace for $\A_n = \Alg
\{ J_+, J_-, J_z\}$.

To see minimality, fix $j,\mu$ and let  $\ket{\psi}$ be a non-zero
vector inside $\H(j,\mu)$. Then by Lemma~\ref{weighted} there is a
$p\geq 0$ such that $J_-^p \ket{\psi}$ is a non-zero multiple of
$\ket{j,-j,\mu}$. Hence, each basis vector $\ket{j,m,\mu}$, for
$-j\leq m \leq j$, belongs to the subspace $\A_n \ket{\psi}
=\H(j,\mu)$, and it follows that $\H(j,\mu)$ is minimal
$\A_n$-reducing. \bx

The structure of $\Delta_n$ determines which of the $\H(j,\mu)$
sum to give the family of minimal central projections. Recall that
the minimal central projections of $\A_n$ and $\A_n^\prime$ are
the same since $\A_n\cap\A_n^\prime =
(\A_n^\prime)^\prime\cap\A_n^\prime$.

\begin{lem}\label{mincentral}
For each $j,\mu$ let $P_{j,\mu}$ be the projection of $\H$ onto
$\H(j,\mu) \equiv P_{j,\mu} \H$. Then the minimal central
projections for $\A_n$ and $\A_n^\prime$ are $\{ P_j\}$ where
\[
P_j = \sum_\mu P_{j,\mu},
\]
and hence $\W_j = P_j \H = \sum_\mu \oplus P_{j,\mu} \H$.
\end{lem}

\Prf The projections $P_{j,\mu}$ form the (unique) maximal family
of non-zero minimal reducing projections for $\A_n$. Thus, the
minimal central projections for $\A_n$ are given by sums of the
$P_{j,\mu}$, and so we must find which subsets of the $P_{j,\mu}$
are `linked inside $\A_n$'. Since linked projections amongst the
$P_{j,\mu}$ necessarily have the same rank, it is enough to fix
$j$ and consider the family $\{P_{j,\mu}\}_\mu$.

In fact, we claim that the entire family $\{P_{j,\mu}\}_\mu$ is
linked inside $\A_n$. To see this, it is sufficient, and best for
use in the noiseless subsystem method, to exhibit bases for
$P_{j,\mu}\H$ which allow us to view the links explicitly. By
construction, the basis $\{\ket{j,m,\mu}: -j\leq m \leq j\}_\mu$
for $P_{j,\mu}\H$ used in the analysis above is such a basis.
Indeed, we may compute that
\begin{eqnarray}\label{link}
\bra{j,m,\mu_1} J_-^{p_1} A J_+^{p_2} \ket{j,m,\mu_1} =
\bra{j,m,\mu_2} J_-^{p_1} A J_+^{p_2} \ket{j,m,\mu_2},
\end{eqnarray}
for all possible choices of $\mu_1,\mu_2,p_1,p_2$ and $A\in\A_n$.
Recall that $\A_n$ is generated by $J_+,J_-,J_z$ as an algebra. By
design, (\ref{link}) is evident for $A$  equal to one of these
generators, for any monomial in them, and hence, when extending by
linearity, for any element of $\A_n$. It follows that for all $j$,
the projection $P_j = \sum_\mu P_{j,\mu}$ is a minimal central
projection for $\A_n$ (and $\A_n^\prime$). \bx

\vspace{0.1in}

{\noindent}{\it Proof of Theorem~\ref{commthm}.} By the previous
result $\A_n$ has a block diagonal decomposition $\A_n =
\sum_{j\in\J_n} \oplus \, \A_{(j)}$, where each $\A_{(j)}$ is a
{\it subalgebra} of $\A_n$ which is unitarily equivalent to
$\A_{(j)}\simeq \one_{p_j}\otimes \M_{q_j}$, since $\rank
P_{j,\mu} = q_j = \dim P_{j,\mu}\H$ for all $\mu$ and there are
$p_j$ linked projections $\{ P_{j,\mu}\}_\mu$. Therefore, the
commutant $\fix(\E_n) = \A_n^\prime$ may be obtained by
\[
\fix(\E_n) = \A^\prime_n = \sum_{j\in\J_n} \oplus \,
\A^\prime_{(j)},
\]
with $\A_{(j)}^\prime \simeq \M_{p_j} \otimes \one_{q_j}$ for
$j\in\J_n$, as claimed in the statement of Theorem~\ref{commthm}.
Observe that we also have the minimal reducing projections for
$\A_n^\prime$ which are supported on the minimal central
projections $P_j$. For each $j\in\J_n$ they are the projections of
rank $p_j$ onto $\spn\{ \ket{j,m,\mu}\}_\mu$. Thus, the explicit
spatial decomposition of $\fix(\E_n) = \A_n^\prime$ is now
evident. \bx

The following is a consequence of the work in this section.

\begin{cor}\label{jtotal}
The set of spectral projections for $J^2$ coincides with the set
of minimal central projections for $\A_n^\prime$ and $\A_n$.
\end{cor}

%
%


\section{Generalized Collective Rotation Channels}\label{S:general}

In this section we consider natural generalizations of collective
rotation channels to higher dimensional representations of $su(2)$
(see Note~\ref{repn}). The commutation relations satisfied by the
Pauli matrices are the defining properties of the Lie algebra
$su(2)$. So far, we have restricted our attention to the special
case where this algebra is represented by $2\times 2$ complex
matrices; specifically the Pauli matrices.  Nevertheless, the
algebra $su(2)$ has an irreducible representation for every
integer dimension; i.e., given $d\geq 1$ it is possible to find
three matrices $\Sigma_{x,d},\ \Sigma_{y,d},\ \Sigma_{z,d}$ of
dimension $d$ satisfying the Pauli commutation relations. Hence,
the rotation group $SU(2)$ also has a representation in every
integer dimension.

Note that the operators $J_x$, $J_y$, $J_z$, which act on
$2^n$-dimensional space, form a representation of the Lie group
$su(2)$. But it is not an irreducible representation as the
theorem in the last section shows; thus the existence of noiseless
subsystems. The irreducible representations of $su(2)$ are
determined by restricting these operators to a minimal reducing
subspace $\H(j,\mu)$. Indeed, it is easily seen that these are
$q_j$-dimensional irreducible representations of the Lie algebra
$su(2)$. (The restrictions of $J_x$, $J_y,$ $J_z$ to each of these
irreducible subspaces satisfies the Pauli commutation relations.)

Physicists call a $d$-dimensional representation of $su(2)$ a
`spin-$s$' representation, where $d=2s+1$.  Hence, the spin $s =
\frac{d-1}{2}$ can take integer and half integer values. From this
more general perspective, we see that in the previous section we
considered the spin-$\frac 12$ ($d=2$) representation of the
rotation group acting on the 2-dimensional Hilbert space of a
qubit. Consideration of the proof in the previous section shows
that it primarily depends on the commutation relations satisfied
by the generators of $su(2)$, not the particular representations
of eigenvectors used in the proof. This `coordinate-free' approach
allows us to readily generalize our results to collective rotation
channels of arbitrary integer dimension. Most of the results from
the previous section follow with small modifications, thus we
shall only outline the approach.

First let us establish some notation. Let $\Sigma_{k,d}$,
$k=x,y,z$, be $(2s+1) \times (2s+1)$ complex matrices forming an
irreducible representation of $su(2)$. These matrices act on the
$d$-dimensional Hilbert space $\H_d$ of a `qudit', where $d=2s+1$.
Consider a collection of $n$ qudits, and their associated
collective rotation generators $J_{x,d},\ J_{y,d},\ J_{z,d}$ on
$\H_d^{\otimes n}$, where, for instance $J_{x,d} = \sum_{k=1}^n
J_{x,d}^{(k)}$ and $J_{x,d}^{(k)} = \ldots \otimes \one_d \otimes
\Sigma_{x,d} \otimes \one_d \ldots$ with $\Sigma_{x,d}$ in the
$k$th tensor slot. As before we may define a unital channel
\[
\E_{n,d}(T) = E_{x,d}TE_{x,d}^\dagger + E_{y,d}TE_{y,d}^\dagger +
E_{z,d}TE_{z,d}^\dagger
\]
where $E_{x,d} = \exp(i\theta_x J_{x,d})$, etc. Let
\[
\A_{n,d} = \Alg\{E_{x,d},E_{y,d},E_{z,d}\} =
\Alg\{J_{x,d},J_{y,d},J_{z,d}\},
\]
the interaction algebra for the channel. Thus the noise commutant
and fixed point set coincide; $\fix(\E_{n,d}) = \A^\prime_{n,d}$.

\begin{prop}
The eigenvalues of $\Sigma_{z,d}$ are $-s, -s+1, \ldots s$, where
$s = \frac{d-1}{2}$.
\end{prop}

\Prf This follows from the definition of $\Sigma_{z,d}$ as the
restriction of $J_z$ on $\H(s,\mu)$. \bx

As in the qubit case ($s=\frac 12$, $d=2$), we can thus represent
a vector in $\H_d^{\otimes n}$ by $\ket{\vec i} =
\ket{i_1i_2,\ldots i_n}$ where $i_k \in \{-s,-s+1,\ldots s\}$
denotes the eigenvalue of $\Sigma_{z,d}$ on the $k$th qudit. With
this notation, we can restate Lemma~\ref{Jzeigen} for arbitrary
finite dimension $d$.

\begin{lem}
For $m = -sn, sn+1,\ldots, sn$ consider the subspaces of
$\H_{d^n}$ given by
\[
\V_m = \spn \big\{ \ket{\vec{i}} : | \vec{i} | =m \big\},
\]
where $|\vec{i}| = \sum_{j=1}^n i_j$. Then $\H_{d^n} = \sum_{-sn
\leq m\leq sn} \oplus \V_m$ and
\[
\dim \V_m = \sum_{k_1 + \ldots + k_n = m+ns} {n \choose k_1\cdots
k_n},
\]
where $k_i\in\{0,\ldots, d-1\}$ and no repeats are allowed, even
reordering, amongst the $n$-tuples $(k_1,\ldots ,k_n)$. Further,
$\V_m$ is an eigenspace for $J_z$ corresponding to the eigenvalue
$m$.
\end{lem}

The proof of this Lemma follows exactly the same lines as
Lemma~\ref{Jzeigen}. The analogues of  Lemmas~\ref{jplusjz} and
\ref{J2eigen} also follow in a straightforward manner; they only
involve the commutation relations which are independent of the
representation of the algebra.

We can thus construct a basis for $\H_{d^n}$ by generalizing the
previous construction. The basis states are $\ket{j,m,\mu}$. The
label $j$ is for the eigenspaces of the operator $(J^{(d)})^2
\equiv J_{x,d}^2 + J_{y,d}^2 + J_{z,d}^2$ which has eigenvalues
given by $j^2 + j$ with $j\in \J_{n,d}$ where
\[
\J_{n,d} = \left\{
\begin{array}{cl}
\{0,1,\ldots, ns \} & \mbox{if $ns$ is an integer} \\
\{\frac{1}{2}, \frac{3}{2}, \ldots, ns \} & \mbox{if $ns$ is a
half integer}
\end{array}\right.
\]
The eigenspaces of $J_{z,d}$ are labelled by $m$, where $m=-j,
-j+1, \ldots, j$ (Recall that $(J^{(d)})^2$ and  $J_{z,d}$
commute, so they can be simultaneously diagonalized.) Finally,
$\mu$ is the extra index required to construct a basis in the
common eigenspace of $(J^{(d)})^2$ and $J_{z,d}$ determined by a
given pair $j,m$.

Let us construct these states as we did in the previous section.
We start with the state $\ket{ns,-ns,1}$ which is the unique
eigenvector of $J_{z,d}$ with eigenvalue $-ns$. It is thus an
eigenvector of $(J^{(d)})^2$.  Then, $J_{+,d} \ket{ns,-ns,1}$ is
an eigenstate of $J_{z,d}$ with eigenvalue $-ns+1$. Furthermore,
since $J_{+,d}$ commutes with  $(J^{(d)})^2$, the vectors $J_{+,d}
\ket{ns,-ns,1}$ and $\ket{ns,-ns,1}$ are in the same  eigenspace
of $(J^{(d)})^2$, hence after normalizing  we can label $J_{+,d}
\ket{ns,-ns,1}$ by $\ket{ns,-ns+1,1}$. By repeating this
procedure, we find an orthonormal basis for the space
\begin{eqnarray*}
\H(ns,1) &=& \spn\{(J_{+,d})^p\ket{ns,-ns,1}: 0\leq p\leq 2ns\} \\
&=& \spn\{\ket{ns,m,1}:m=-ns,-ns+1,\ldots,ns\}.
\end{eqnarray*}
By construction, $\H(ns,1)$ is a minimal reducing subspace for
$\A_{n,d}$.  Furthermore, since the spectral projections of
$(J^{(d)})^2$ are the minimal central projectors of $\A_{n,d}$,
the subspace $\H(ns,1)$ is an eigenspace of $(J^{(d)})^2$.

We then consider the subspace $\V_{-ns+1} \ominus \spn
\{\ket{ns,-ns+1,1}\}$. This is the eigenspace of $J_{z,d}$ with
eigenvalue $m=-ns+1$ which is perpendicular to the eigenspace of
$(J^{(d)})^2$ labelled by $ns$. Hence, these vectors require a
different $j$ label, say $j=ns-1$. We can now choose a basis for
$\V_{-ns+1} \ominus \spn\{\ket{ns,-ns+a,1}\}$, which is labeled
$\ket{ns-1,-ns+1,\mu}$ where the first two terms just label the
subspace $\V_{-ns+1} \ominus \spn\{\ket{ns,-ns+1,1}\}$ and $\mu$
is an extra label to form a basis within this subspace. Thus, as
we did in the previous section, we construct subspaces by applying
the shift operator
\begin{eqnarray*}
\H(ns-1,\mu) &=& \spn\{J_{+,d}^p\ket{ns-1,-ns+1,\mu}: 0\leq p\leq 2(ns-1)\} \\
&=& \spn\{\ket{ns-1,m,\mu}:-ns+1 \leq m \leq ns-1\}.
\end{eqnarray*}

This procedure can be repeated with the subspaces
\[
\V_{-m} \ominus \spn \{\ket{j,-m,\mu}:j=m+1,\ldots ns, \mu =
1,\ldots, q_j\}
\]
to form the subspaces
\begin{eqnarray*}
\H(j,\mu) &=& \spn\{J_{+,d}^p\ket{j,-j,\mu}: 0\leq p\leq 2j\} \\
&=& \spn\{\ket{j,m,\mu}:-j \leq m \leq j\}.
\end{eqnarray*}
The subspaces $\H(j,\mu)$ are minimal $\A_{n,d}$-reducing and for
fixed $j$, the subspaces $\{ \H(j,\mu)\}_{\mu}$ are linked inside
$\A_{n,d}$. Thus with this analysis in hand, we may state the
following generalization of Theorem~\ref{commthm}.

\begin{thm}\label{commthm1}
Let $\E_{n,d}$ be the collective rotation channel for  fixed
positive integers $n\geq 1$ and $d\geq 2$. Then
\begin{eqnarray}\label{collfixed}
\fix (\E_{n,d}) = \A_{n,d}^\prime = \sum_{j\in\J_{n,d}} \oplus \,
\A_{(j)}^\prime,
\end{eqnarray}
where $\A_{(j)}^\prime$ is a $\ca$-subalgebra of $\A_{n,d}^\prime$
given, up to unitary equivalence, by
\[
\A_{(j)}^\prime \simeq \M_{p_j} \otimes \one_{q_j} \qfor
j\in\J_{n,d},
\]
with $p_{ns} = 1$ where and for $j\in\J_{n,d}$, $j< ns$,
\[
p_j = \dim \V_j - \dim \V_{j+1}
\]
where
\[
q_j = 2j +1 \qfor j\in\J_{n,d}.
\]
\end{thm}

\begin{rem}
In  light of this analysis, we can extend the result to more
general Lie groups. Let $\G$ be a compact connected semisimple Lie
group and $\G^{\otimes n}$ denote its $n$-fold tensor product.
Further, let $\Sigma_k$ be the set of generators of the associated
Lie algebra. This algebra is entirely specified by its {\em
structure constants} $C_{kmn}$ defined by
\begin{equation}
[\Sigma_m,\Sigma_n] = i\sum_k C_{kmn}\Sigma_k.
\end{equation}
The operators
\begin{equation}
J_k = (\Sigma_k\otimes\one \otimes\one\otimes\ldots) +
(\one\otimes\Sigma_k\otimes\one\otimes\ldots) + \ldots
\end{equation}
are generators of the generalized `collective rotation' which is a
subgroup of $\G^{\otimes n}$. Clearly, they have the same
structure constants has the $\Sigma_k$; they represent the same
algebra. Nevertheless, the $J_k$ do not form an irreducible
representation of the algebra. Hence, it is possible to write them
as a direct sum of irreducible representations. A special property
of these representations is that all the projections onto the
irreducible subspaces of the same dimension are in fact `linked'
inside the algebra.  Thus, it follows that there is an abundance
of noiseless subsystems which can be explicitly identified for the
corresponding quantum channels. An expansion of this analysis is
contained in \cite{JKK}.
\end{rem}


{\noindent}{\it Acknowledgements.}  We would like to thank Marius
Junge and Peter Kim  for helpful  conversations. We are grateful
for support from the Perimeter Institute, the Institute for
Quantum Computing, the University of Guelph, NSERC, MITACS, and
ARDA.





\end{document}